\documentclass[12pt,a4paper]{article}

\unitlength\textwidth
\divide\unitlength by 300\relax

\usepackage{amsmath,amssymb}
\usepackage{bm}
\bmdefine{\NNN}{N}
\bmdefine{\ZZZ}{Z}
\bmdefine{\RRR}{R}
\bmdefine{\aaa}{a}
\bmdefine{\bbb}{b}
\bmdefine{\eee}{e}
\bmdefine{\zerovec}{0}
\bmdefine{\onevec}{1}
\bmdefine{\alphaaa}{\alpha}
\bmdefine{\rhooo}{\rho}

\newcommand{\mmmm}{\mathfrak{m}}

\newcommand{\qed}{\nolinebreak\rule{.3em}{.6em}}

\newcommand{\supp}{\mathrm{supp}}
\newcommand{\Tot}{\mathrm{Tot}}

\newtheorem{thm}{Theorem}[section]

\newtheorem{remark}[thm]{Remark}

\makeatletter
\newcommand{\mysloppy}{\tolerance 9999 \hfuzz .5\p@ \vfuzz .5\p@}
\makeatother
\begin{document}\mysloppy
\begin{center}
\LARGE
Behavior of local cohomology modules under polarization
\end{center}

\begin{center}
\large
Mitsuhiro \sc Miyazaki%
\footnote{%
Dept.\ Math.\
Kyoto University of Education,
Fushimi-ku, Kyoto, 612-8522 Japan
\ \ 
E-mail:
\tt
g53448@kyokyo-u.ac.jp}
\end{center}

\begin{trivlist}\item[]\small
{\bf Abstract:}\quad
Let $S=k[x_1,\ldots, x_n]$ be
a polynomial ring over a field $k$ with $n$ variables
$x_1$, \ldots, $x_n$,
$\mmmm$ the irrelevant maximal ideal of $S$,
$I$ a monomial ideal in $S$ and
$I'$ the polarization of $I$ in the polynomial ring
$S'$ with $\rho$ variables.
We show that each graded piece 
$H_\mmmm^i(S/I)_\aaa$, $\aaa\in\ZZZ^n$,
of the local cohomology module $H_\mmmm^i(S/I)$
is isomorphic to  a specific graded piece
$H_{\mmmm'}^{i+\rho-n}(S'/I')_\alphaaa$,
$\alphaaa\in\ZZZ^\rho$, of the local cohomology module
$H_{\mmmm'}^{i+\rho-n}(S'/I')$,
where $\mmmm'$ is the irrelevant maximal ideal of $S'$.
\end{trivlist}
\begin{trivlist}\item[]\small
\end{trivlist}

\section{Introduction}

Let $S=k[x_1,\ldots, x_n]$ be
a polynomial ring over a field $k$ with $n$ variables
$x_1$, \ldots, $x_n$
and $\mmmm$ the irrelevant maximal ideal of $S$.
For a monomial ideal $I$ in $S$, the local cohomology modules
$H_\mmmm^i(S/I)$ have $\ZZZ^n$-graded structure.

Hochster described each graded piece of $H_\mmmm^i(S/I)$
by using the reduced cohomology group of a simplicial complex
related to $I$ when $I$ is  square-free (see \cite[II 4.1 Thoerem]{sta}).
Takayama \cite[Theorem 1]{tak} generalized this result to the
case where $I$ is not necessarily square-free.

On the other hand, there is a technique 
which associates a not necessarily square-free monomial ideal $I$
with a square-free monomial ideal, 
called the polarization of $I$,
sharing many ring theoretical properties with $I$.

In this note, we show that the each graded piece of the
local cohomology modules $H_\mmmm^i(S/I)$
is isomorphic to a specific graded piece of the local cohomology
module of the polarization of $I$.

\section{Preliminaries}
Let $S=k[x_1,\ldots, x_n]$ be
a polynomial ring over a field $k$ with $n$ variables
$x_1$, \ldots, $x_n$
and $\mmmm$ the irrelevant maximal ideal of $S$
and $I$ a monomial ideal of $S$.

For a monomial $m=x_1^{b_1}\cdots x_n^{b_n}$ in $S$,
we set $\nu_i(m)=b_i$.
We denote by $G(I)$ the minimal set of monomial generators of $I$.
Set $\rho_i=\max\{\nu_i(m)\mid m\in G(I)\}$
for $i=1$, $2$, \ldots, $n$
and $\rho=\rho_1+\cdots+\rho_n$.
Then the polarization of $I$ is defined as follows.

Let $S'=k[x_{ij}\mid 1\leq i\leq n$, $1\leq j\leq \rho_i]$
be the polynomial ring with $\rho$ variables $\{x_{ij}\}$.
For a monomial $m$ in $S$, we set $m'=\prod_{i=1}^n\prod_{j=1}^{\nu_i(m)}x_{ij}$.
Then the polarization $I'$ of $I$ is defined by
$I'=(m'\mid m\in G(I))S'$.
It is clear from the definition that $I'$ is a square-free
monomial ideal.
Furthermore,
it is known that 
$
\{x_{ij}-x_{i1}\mid 1\leq i\leq n, 2\leq j \leq \rho_i\}
$
is an $S'/I'$-regular sequence in any order and
\[
S'/(I'+(x_{ij}-x_{i1}\mid 1\leq i\leq n, 2\leq j \leq \rho_i))\simeq S/I.
\]

For vectors $\aaa=(a_1,\ldots, a_n)$ and $\bbb=(b_1,\ldots, b_n)$,
we denote $\aaa\leq\bbb$ to express that $a_i\leq b_i$ for $i=1$, \ldots, $n$.
And we define $\supp_-\aaa=\{i\mid a_i<0\}$
and call the negative support of $\aaa$.
We set $\zerovec=(0,0,\ldots,0)$, $\onevec=(1,1,\ldots,1)$
and $\rhooo=(\rho_1,\rho_2,\ldots, \rho_n)$.

We denote the cardinality of a finite set $X$ by $|X|$
and the set $\{1,2, \ldots, n\}$ by $[n]$.

Here we recall the result of Takayama.

\begin{thm}[{\cite[Theorem 1]{tak}}]
\label{tak thm}
Let $S$ and $I$ be as above and $\aaa\in\ZZZ^n$.
Set
$\Delta_\aaa=\{F\setminus\supp_-\aaa\mid
[n]\supset F\supset\supp_-\aaa,
\forall m\in G(I)\exists i\in[n]\setminus F;a_i<\nu_i(m)\}$.
Then 
\[
H_\mmmm^i(S/I)_\aaa\simeq\tilde H^{i-|\supp_-\aaa|-1}(\Delta_\aaa;k).
\]
\end{thm}
Note that $\Delta_\aaa$ is a simplicial complex with vertex set $[n]$.
We call $\Delta_\aaa$ the Takayama complex.
Note also that if $a_i\geq \rho_i$, then $\Delta_\aaa$ is
a cone over $i$ and all the reduced cohomology modules vanish.
Therefore
$H_\mmmm^i(S/I)_\aaa=0$ if $\aaa\not\leq\rhooo-\onevec$.

\section{Main theorem}

Now we state the main result of this paper.

\begin{thm}\label{main thm}
With the notation in previous section,
assume that $\aaa\leq\rhooo-\onevec$.
Set
\[
\alpha_i=
\left\{
\begin{array}{ll}
(\underbrace{0,\ldots,0}_{a_i+1},\underbrace{-1,\ldots,-1}_{\rho_i-a_i-1})&
\quad\mbox{if $a_i\geq 0$,}\\
(\underbrace{-1,\ldots,-1}_{\rho_i})&
\quad\mbox{if $a_i< 0$}
\end{array}
\right.
\]
and
$\alphaaa=(\alpha_1,\alpha_2,\ldots,\alpha_n)\in\ZZZ^\rho$.
Then
\[
H_\mmmm^i(S/I)_\aaa\simeq
H_{\mmmm'}^{i+\rho-n}(S'/I')_\alphaaa,
\]
where $\mmmm'$ is the irrelevant maximal ideal of $S'$.
\end{thm}
{\bf proof.}\ \ 
{\bf Step 1.}
We first consider the case where $\aaa=\rhooo-\onevec$.

First note that local cohomology modules $H_{\mmmm'}^i(S'/I')$
have not only the $\ZZZ^\rho$-grading but also the
$\ZZZ^n$-grading by setting $\deg x_{ij}=\eee_i$,
where $\eee_i$ is the $i$-th fundamental vector in $\ZZZ^n$.

We denote by $C^\bullet$ the \v{C}ech complex with respect to
$\{x_{ij}\mid 1\leq i\leq n$, $1\leq j\leq \rho_i\}$
and $K_\bullet$ the Koszul complex with respect to
$\{x_{ij}-x_{i1}\mid 1\leq i\leq n, 2\leq j \leq \rho_i\}$.
Set $M^{p,q}=C^p\otimes_{S'}K_{\rho-n-q}\otimes_{S'}S'/I'$.
Then $M^{\bullet,\bullet}$ is a third quadrant double complex
which has a $\ZZZ^n$-graded structure.
Set $\{'E_r\}$ and $\{''E_r\}$ be the spectral sequences arising
from $M^{\bullet,\bullet}$.

Since
$\{x_{ij}-x_{i1}\mid 1\leq i\leq n, 2\leq j \leq \rho_i\}$
is an $S'/I'$-regular sequence, we see that
\[
'E_1^{p,q}
\simeq
\left\{
\begin{array}{ll}
C^p\otimes S/I&\quad\mbox{$q=\rho-n$,}\\
0&\quad\mbox{otherwise.}
\end{array}
\right.
\]
And the horizontal complex $'E_1^{\bullet,\rho-n}$ is isomorphic to
the \v{C}ech complex with respect to $\underbrace{x_1, x_1, \ldots, x_1}_{\rho_1}$, 
$\underbrace{x_2, x_2, \ldots, x_2}_{\rho_2}$, \ldots, 
$\underbrace{x_n, x_n, \ldots, x_n}_{\rho_n}$.
Therefore
\[
'E_2^{p,q}
\simeq
\left\{
\begin{array}{ll}
H_\mmmm^p(S/I)\quad&\mbox{$q=\rho-n$,}\\
0&\mbox{otherwise.}
\end{array}
\right.
\]
So the spectral sequence $\{'E_r\}$ collapses and we see that
\[
H^i(\Tot(M^{\bullet,\bullet}))\simeq 
{}'E_2^{i-\rho+n,\rho-n}\simeq
H_\mmmm^{i-\rho+n}(S/I)
\]
for any $i\in\ZZZ$,
where $\Tot(M^{\bullet,\bullet})$ is the total complex of $M^{\bullet,\bullet}$.

Next we consider $\{''E_r\}$.
It is clear that
$''E_1^{p,q}\simeq K_{\rho-n-q}\otimes_{S'}H_{\mmmm'}^p(S'/I')$.
Since $I'$ is square-free, $H_{\mmmm'}^p(S'/I')_{\alphaaa}=0$
 if $\alphaaa\not\leq\zerovec$
by the remark after Theorem \ref{tak thm}.
Therefore
$H_{\mmmm'}^p(S'/I')_l=0$ if $l>0$, $l\in\ZZZ$,
where $H_{\mmmm'}^p(S'/I')_l$ denotes the total degree 
$l$ piece of $H_{\mmmm'}^p(S'/I')$.

Since $K_{\rho-n-q}$ is a free $S'$-module with free basis 
consisting of 
total degree $\rho-n-q$ elements,
we see that
\[
(''E_1^{p,q})_{\rho-n}=0\quad\mbox{if $q\neq0$.}
\]
Therefore total degree $\rho-n$ piece of $\{''E_r\}$
collapses and we see that
\[
H^i(\Tot(M^{\bullet,\bullet}))_{\rho-n}\simeq
(''E_1^{i,0})_{\rho-n}\simeq
H_{\mmmm'}^i(S'/I')_0
\]
for any $i\in\ZZZ$.

So 
\[
H_{\mmmm}^{i-\rho+n}(S/I)_{\rho-n}\simeq
H_{\mmmm'}^i(S'/I')_0
\]
for any $i\in\ZZZ$.
By the remark after Theorem \ref{tak thm},
we see that 
$H_\mmmm^j(S/I)_{\rho-n}=H_\mmmm^j(S/I)_{\rhooo-\onevec}$
and
$H_{\mmmm'}^j(S'/I')_0=H_{\mmmm'}^j(S'/I')_\zerovec$
for any $j\in\ZZZ$.
This means 
\[
H_\mmmm^i(S/I)_{\rhooo-\onevec}\simeq
H_{\mmmm'}^{i+\rho-n}(S'/I')_\zerovec
\]
and
this is what we wanted to prove.

{\bf Step 2.}
Next we consider the case where
$a_i=\rho_i-1$ or $a_i<0$ for any $i=1$, \ldots, $n$.

By changing the subscripts, we may assume that
$a_i=\rho_i-1$ for $i=1$, \ldots, $m$ and 
$a_i<0$ for $i=m+1$, \ldots, $n$.
Set
$S_+=k[x_1,\ldots, x_m]$,
$I_+=IS[x_{m+1}^{-1},\ldots, x_n^{-1}]\cap S_+$
and
$\mmmm_+$ the irrelevant maximal ideal of $S_+$.
Note that $I_+$ is the monomial ideal of $S_+$ generated by
the monomials obtained by substituting 1 to $x_{m+1}$, \ldots, $x_n$
of the monomials in $I$.
Note also that if we denote the Takayama complex with respect to
$I_+$ and $(a_1,\ldots, a_m)$ by $\Delta_+$,
then $\Delta_+=\Delta_\aaa$.
Therefore
\[
\begin{array}{rcl}
H_\mmmm^i(S/I)_\aaa
&\simeq&\tilde H^{i-|\supp_-\aaa|-1}(\Delta_\aaa;k)\\
&\simeq&\tilde H^{i-(n-m)-1}(\Delta_+;k)\\
&\simeq&H_{\mmmm_+}^{i-n+m}(S_+/I_+)_{(a_1,\ldots,a_m)}
\end{array}
\]
since 
$|\supp_-\aaa|=n-m$ and 
$|\supp_-(a_1,\ldots,a_m)|=0$.

We also set
$S'_+=k[x_{ij}\mid 1\leq i\leq m, 1\leq j\leq \rho_i]$,
$I'_+=I'S'[x_{ij}^{-1}\mid m+1\leq i\leq n, 1\leq j\leq \rho_i]\cap S'_+$
and $\mmmm'_+$ the irrelevant maximal ideal of $S'_+$.
Then it is easily verified that $I'_+$ is the polarization of $I_+$.
And we see that 
$H_{\mmmm'}^{i+\rho-n}(S'/I')_\alphaaa
\simeq H_{\mmmm'_+}^{i+\rho_1+\cdots+\rho_m-n}(S'_+/I'_+)_{(\alpha_1,\ldots,\alpha_m)}$
by the same argument above.
Therefore, we can reduce this case to the case of step 1.

{\bf Step 3.}
Finally, we consider the general case.

Assume that $0\leq a_n<\rho_n-1$.
Consider the ``partial polarization with respect to $x_n$''.
That is, set
$S''=S[x_{nj}\mid a_n+2\leq j\leq \rho_n]$
and for a monomial $m$ in $S$, set
\[
m''=
\left\{
\begin{array}{ll}
m&\mbox{if $\nu_n(m)\leq a_n+1$,}\\
\prod_{i=1}^{n-1}x_i^{\nu_i(m)}x_n^{a_n+1}\prod_{j=a_n+2}^{\nu_n(m)}x_{nj}
\quad&\mbox{if $\nu_n(m)\geq a_n+2$.}
\end{array}
\right.
\]
Set also $I''=(m''\mid m\in G(I))$.
Then the polarization of $I''$ is the same as that of $I$.

And if we denote the Takayama complex with respect to $I''$ and
$(a_1,\ldots, a_{n-1}, a_n, \underbrace{-1,\ldots,-1}_{\rho_n-1-a_n})$
by $\Delta''$,
then we see that $\Delta''=\Delta_\aaa$.
Therefore
\[
\begin{array}{rcl}
H_\mmmm^i(S/I)_\aaa
&\simeq&\tilde H^{i-|\supp_-\aaa|-1}(\Delta_a;k)\\
&\simeq&\tilde H^{i-|\supp_-(a_1,\ldots, a_{n-1}, a_n, -1,\ldots,-1)|-1+(\rho_n-1-a_n)}(\Delta'';k)\\
&\simeq&H_{\mmmm''}^{i+\rho_n-(a_n+1)}(S''/I'')_
{(a_1,\ldots, a_{n-1}, a_n, -1,\ldots,-1)}.
\end{array}
\]
So we can reduce the proof to the case where $a_n=\rho_n-1$.

Using this argument repeatedly, we may assume that
$a_i\geq 0\Rightarrow a_i=\rho_i-1$,
i.e.,
we can reduce the proof of the theorem to the case of
step 2.
\qed

\begin{remark}\rm
Step 1 of the proof of Theorem \ref{main thm}
can also be proved by using \cite[Corollary 5.2]{sba},
insted of the spectral sequence argument.
\end{remark}

\end{document}